\documentclass[10 pt,leqno]{article}
\usepackage{amsfonts,amssymb,amsmath}

\begin{document}

\author{Ajai Choudhry}
\title{A diophantine problem on biquadrates revisited}
\date{}
 \maketitle

\begin{abstract} In this paper we obtain a new parametric solution of the problem of finding two triads of biquadrates with equal sums and equal products.
 
\end{abstract}



\bigskip
\bigskip

This paper is a sequel to my earlier paper \cite{Cho2} concerning the problem of finding two triads of biquadrates with equal sums and equal products. The problem is mentioned in \cite[p. 214]{Gu}, and  thirty numerical solutions, obtained by computer search, were presented in \cite{Cho2}. Subsequently, a method of obtaining infinitely many numerical examples of such triads was given in \cite[pp. 346--47]{ChW}. In this paper we obtain a parametric solution of the problem.

To obtain the desired triads of biquadrates, we will first obtain a parametric solution of  the simultaneous diophantine equations, 
 \begin{align}
 x_1^2+x_2^2+x_3^2&=y_1^2+y_2^2+y_3^2, \label{sys1}\\
x_1x_2x_3&=y_1y_2y_3, \label{sys2}
\end{align}
and then choose the parameters such  that $x_i,\;y_i,\;i=1,\,2,\,3$, become perfect squares. The complete solution of Eqs.~\eqref{sys1} and \eqref{sys2} is already known (\cite{Cho2}, \cite{Ke}).   Since we need to ensure that $x_i,\;y_i,\;i=1,\,2,\,3$ eventually become perfect squares, instead of using the solutions given in (\cite{Cho2}, \cite{Ke}), we  will  solve Eqs.~\eqref{sys1} and \eqref{sys2} directly by writing,
\begin{equation}
\begin{aligned}
x_1&=-(t u+1)p^2,\quad&x_2&=(u+1)q^2,\quad &x_3&=r^2,\\
y_1&=(u+1)p^2,\quad &y_2&=q^2,\quad &y_3&=-(t u +1)r^2,
\end{aligned}
\label{subsxy}
\end{equation}
where $p, \, q,\,r, \,t$ and $u$ are arbitrary nonzero parameters. 

With the values of $x_i,\;y_i,\;i=1,\,2,\,3$ given by \eqref{subsxy}, Eq.~\eqref{sys2} is identically satisfied while Eq.~\eqref{sys1} reduces to a linear equation in $u$ which is readily solved, and we thus obtain a parametric solution of the simultaneous equations \eqref{sys1} and \eqref{sys2}. Since both these equations are homogeneous in the variables $x_i,\;y_i$, we obtain the following solution of Eqs.~\eqref{sys1} and \eqref{sys2} after appropriate scaling:
\begin{equation}
\begin{aligned}
x_1&=\{(p^4-r^4)t^2-2(p^4-q^4)t+p^4-q^4\}p^2k,\\
x_2&=\{(p^4-r^4)t^2-2(p^4-r^4)t+p^4-q^4\}q^2k,\\
x_3&=\{(p^4-r^4)t^2-(p^4-q^4)\}r^2k,\\
y_1&=\{(p^4-r^4)t^2-2(p^4-r^4)t+p^4-q^4\}p^2k,\\
y_2&=\{(p^4-r^4)t^2-(p^4-q^4)\}q^2k,\\
y_3&=\{(p^4-r^4)t^2-2(p^4-q^4)t+p^4-q^4\}r^2k,
\end{aligned} \label{solsys12}
\end{equation}
where $p,\, q,\,r,\,t$ and $k$ are arbitrary parameters.

To ensure that $x_i, \;y_i,\,i=1,\,2,\,3,$ are all perfect squares, we note that since the ratios $y_1/x_2,\,y_2/x_3$ and $y_3/x_1$ are perfect squares, it suffices to ensure that $x_i, \;\,i=1,\,2,\,3,$ are all perfect squares. We will first make the products $x_1x_3$ and $x_2x_3$ simultaneously perfect squares by choice of the parameters $p,\,q,\,r$ and then choose $k$ suitably such that $x_1,\,x_2$ and $x_3$ become perfect squares simultaneously. 

 Now $x_1x_3$ and $x_2x_3$ may be considered as two quartic polynomials in $t$ with coefficients in terms of the parameters $p, \,q,\,r$, and we need to choose a suitable value of $t$ such that both of these quartic polynomials become perfect squares. The problem of making two  quartic polynomials  perfect squares simultaneously is a very difficult diophantine problem and there is no general method of finding a solution. Here we will solve the problem by finding different values of $t$ that make each quartic a perfect square separately, and we will equate two such values of $t$ to obtain a condition that must be satisfied by the parameters $p,\,q,\,r$.  Finally, to obtain a solution, we will find suitable values of $p,\,q,\,r$ that satisfy the condition just obtained. 

We note that the coefficient of $t^4$ in both of our quartic  polynomials are perfect squares. Accordingly, we may divide each quartic polynomial by the respective coefficient of $t^4$  to obtain two monic quartic polynomials  both of which must be made perfect squares simultaneously. Thus, we must find a common value of $t$ that satisfies both of  the following two  diophantine equations for suitably chosen values of $p,,q$ and $r$:
\begin{align}
y_1^2&=t^4-2(p^4-q^4)t^3/(p^4-r^4)+2(p^4-q^4)^2t/(p^4-r^4)^2 \nonumber\\
& \quad \quad \quad -(p^4-q^4)^2/(p^4-r^4)^2, \label{quart1}\\
y_2^2&=t^4-2t^3+2(p^4-q^4)t/(p^4-r^4)-(p^4-q^4)^2/(p^4-r^4)^2. \label{quart2}
\end{align}

A method of obtaining  values of $t$ that make a monic quartic polynomial of $t$ a perfect square has been given by Fermat \cite[p. 639]{Di}. This method, however, leads to very cumbersome solutions of \eqref{quart1} and \eqref{quart2}. To obtain relatively simpler solutions of \eqref{quart1} and \eqref{quart2}, we will apply a theorem given in \cite[Theorem 4.1, p. 789-90]{Cho3} by  which we can  use two known solutions of a quartic equation of  type \eqref{quart1} to obtain a new solution. 

It is readily verified that two solutions of Eq.~\eqref{quart1} are given by
\[(y_1,\,t)=\left(1,\,(q^4-r^4)/(p^4-r^4)\right)\]
\[\mbox{\rm and}\;\; 
(y_1,\,t)=\left((p^2+q^2)/(p^2+r^2),\,2pq(p^2+q^2)(q^2-r^2)/\{(p^2-r^2)(p^2+r^2)^2\}\right),\]
and on applying the aforementioned  theorem, we get the following value of $t$ that provides a solution of Eq.~\eqref{quart1}:
\begin{equation}
\frac{(p^2+q^2)(p^4+2p^3q-2p^2r^2-2pq^3-4pqr^2+q^4+2q^2r^2+2r^4)}{(p^2-r^2)(p^4+2p^3q-2p^2r^2+2pq^3+q^4-2q^2r^2-2r^4)}. \label{valt1}  
\end{equation}

Similarly, two solutions of Eq.~\eqref{quart2} are given by 
\[(y_2,\,t)=\left((p^4-q^4)/(p^4-r^4),\;-(p^4-q^4)(q^4-r^4)/(p^4-r^4)^2,\right)\]
\[\mbox{\rm and}\;\; 
(y_2,\,t)=\left((p^2-q^2)/(p^2-r^2),\;2pr(p^2-q^2)(q^2-r^2)/\{(p^2-r^2)^2(p^2+r^2)\}\right),\]
and using these solutions, we get the following value of $t$ that provides a solution of Eq.~\eqref{quart2}:
\begin{equation}\frac{(p^2+q^2)(p^4-2p^3r-2p^2q^2+4pq^2r+2pr^3+2q^4+2q^2r^2+r^4)}{(p^2-r^2)(p^4-2p^3r-2p^2q^2-2pr^3-2q^4-2q^2r^2+r^4)}.
\label{valt2}
\end{equation}

On equating the two values of $t$ given by \eqref{valt1} and \eqref{valt2}, we get the following condition:
\begin{multline}
(q+r)(p^2+q^2)(q^2+r^2)(p^2+pq-pr+q^2-qr+r^2)\\
\times (p^3-pq^2+pqr-pr^2+q^3-r^3)=0.\label{condt}
\end{multline}

On equating to 0 any of the factors on the left-hand side of Eq.~\eqref{condt}, except the last,  we get a trivial result. When, however, we equate to 0 the last factor, we get a cubic equation in $p,\,q,\,r$ which represents a curve of genus 0 in the projective plane with a singularity at $(-1,\,-1,\,1)$. On drawing an arbitrary straight line through the singular point and taking its intersection  with the cubic curve, we get a rational point which immediately yields the following solution for $p,\,q,\,r$:
\begin{equation}
p = 2a^3-3a^2b+3ab^2-b^3,\;\; q = -a^3+3a^2b-2ab^2+b^3,\;\; r = a^3-ab^2+b^3, \label{valpqr}
\end{equation}   
where $a$ and $b$ are arbitrary parameters. With these values of $p,\,q,\,r$, the two values of $t$ given by \eqref{valt1} and \eqref{valt2} coincide and are given by
\begin{equation} -(5a^4-8a^3b+8a^2b^2-4ab^3+b^4)/\{a(a-b)(5a^2-5ab+2b^2)\}. \label{valt} \end{equation}

On substituting the valus of $p,\,q,\,r$ given by \eqref{valpqr} and the value of $t$ given by \eqref{valt} in \eqref{solsys12} and choosing $k$ suitably, we  obtain the following solution of Eqs.~\eqref{sys1} and \eqref{sys2}:
\begin{equation}
\begin{aligned}
x_1& = (2a-b)^2(a^2-ab+b^2)^2(5a^3-7a^2b+4ab^2-b^3)^2, \\
x_2& = (a^3-3a^2b+2ab^2-b^3)^2(5a^3-8a^2b+5ab^2-b^3)^2,\\
 x_3& = (a^3-ab^2+b^3)^2(3a^2-3ab+b^2)^2b^2, \\
y_1& = (2a-b)^2(a^2-ab+b^2)^2(5a^3-8a^2b+5ab^2-b^3)^2, \\
y_2& = (a^3-3a^2b+2ab^2-b^3)^2(3a^2-3ab+b^2)^2b^2, \\
y_3& = (a^3-ab^2+b^3)^2(5a^3-7a^2b+4ab^2-b^3)^2,
\end{aligned} \label{solxy}
\end{equation}
where $a$ and $b$ are arbitrary  parameters.

Since the values of $x_i,\;y_i,\;i=1,\,2,\,3$, given by \eqref{solxy} are all perfect squares, each of the two triads $\{x_1^2,\,x_2^2,\,x_3^2\}$ and $\{y_1^2,\,y_2^2,\,y_3^2\}$ consists of three biquadrates, and in view of the relations \eqref{sys1} and \eqref{sys2}, the two triads of biquadrates have equal sums and equal products.

As a numerical example, when $a=1,\,b=-1$, we get the two triads of biquadrates $\{7^4,\,133^4,\,153^4\}$ and $\{17^4,\,49^4,\,171^4\}$ with equal sums and equal products.

\begin{center}
\Large
Acknowledgment
\end{center}
 
I wish to  thank the Harish-Chandra Research Institute, Allahabad for providing me with all necessary facilities that have helped me to pursue my research work in mathematics.

\bigskip

\noindent Postal Address: 13/4 A Clay Square, \\
\hspace{1.1in} Lucknow - 226001,\\
\hspace{1.1in}  INDIA \\
\medskip
\noindent e-mail address: ajaic203@yahoo.com

\end{document}